\documentclass{amsart}
\usepackage{amssymb,amsmath,latexsym}
\usepackage{amsthm}
\usepackage{fontenc}
\usepackage{amssymb}
\numberwithin{equation}{section}

 \DeclareMathOperator{\spt}{spt}

\begin{document}
\author{Alexander E Patkowski}
\title{Some implications of the $2$-fold Bailey lemma}
\date{\vspace{-5ex}}
\maketitle
\begin{abstract}The $2$-fold Bailey lemma is a special case of the $s$-fold Bailey lemma
introduced by Andrews in 2000.
We examine this special case and its applications to partitions and
recently discovered $q$-series identities.
Our work provides a general comparison of the utility of the $2$-fold Bailey
lemma and the more widely applied $1$-fold Bailey lemma.
We also offer a discussion of the $\spt_M(n)$ function and related identities.\end{abstract}


\rm
\section{Introduction}
The Symmetric Bilateral Bailey transform [4], states that if
\begin{equation}B_n=\sum_{j=-n}^{n}A_ju_{n-j}v_{n+j},\end{equation}
and
\begin{equation}\gamma_n=\sum_{j=|n|}^{\infty}\delta_ju_{j-n}v_{n+j},\end{equation}
then
\begin{equation}\sum_{n=-\infty}^{\infty}A_n\gamma_n=\sum_{n=0}^{\infty}B_n\delta_n.\end{equation}
Here we say that $(A_n,B_n)$ is a Bailey pair (in the symmetric sense), and $(\gamma_n, \delta_n)$ is a conjugate Bailey pair. If we break symmetry, we say a pair
$(\alpha_n, \beta_n)$ is a Bailey pair if
\begin{equation}\beta_n(a)=\sum_{j=0}^{n}\alpha_ju_{n-j}v_{n+j}.\end{equation}
Bailey [6] chose
 $v_n=(aq)_n^{-1},$ $u_n=(q)_n^{-1},$ (where we use standard notation [8] for $q$-shifted factorials) giving the conjugate pair $(\gamma_n, \delta_n)$
 \begin{equation} \delta_n=(\rho_1)_n(\rho_2)_n(aq/\rho_1\rho_2)^n,\end{equation}
 \begin{equation} \gamma_n=\frac{(aq/\rho_1)_{\infty}(aq/\rho_2)_{\infty}}{(aq)_{\infty}(aq/\rho_1\rho_2)_{\infty}}\frac{(\rho_1)_n(\rho_2)_n(aq/\rho_1\rho_2)^n}{(aq/\rho_1)_n(aq/\rho_2)_n}.\end{equation}
 \\*
 In [3] Andrews considered an $s$-fold extension of the Bailey lemma, and gave the definition of the $s$-fold Bailey pair $(A_{n_1,n_2,\dots,n_s},B_{n_1,n_2,\dots,n_s})$ relative to $(a_1,a_2,\dots,a_s),$
 $$B_{n_1, \dots,n_s}=\sum_{r_1=-n_1}^{n_1}\cdots\sum_{r_s=-n_s}^{n_s} \frac{A_{r_1,\dots,r_s}}{(a_1q)_{n_1+r_1} (q)_{n_1-r_1}\dots(a_sq)_{n_s+r_s} (q)_{n_s-r_s}}.$$
Here we focus on the $s=2$ case, in which case Andrews' generalisation of (1.4) is given by
$$B_{n_1, n_2}=\sum_{r_1=-n_1}^{n_1}\sum_{r_2=-n_2}^{n_2} \frac{A_{r_1,r_2}}{(a_1q)_{n_1+r_1} (q)_{n_1-r_1}(a_2q)_{n_2+r_2} (q)_{n_2-r_2}}.$$
For this symmetric $2$-fold Bailey pair relative to $(a_1,a_2)$ we have the identity
\begin{equation}\sum_{n_1=-\infty}^{\infty}\sum_{n_2=-\infty}^{\infty} A_{n_1, n_2}\gamma_{n_1}(a_1)\gamma_{n_2}(a_2)=\sum_{n_1=0}^{\infty}\sum_{n_2=0}^{\infty}B_{n_1,n_2}\delta_{n_1}(a_1)\delta_{n_2}(a_2),\end{equation}
where both $(\gamma_{n_1}(a_1),\delta_{n_1}(a_1))$ and $(\gamma_{n_2}(a_2),\delta_{n_2}(a_2))$ are given by (1.5)--(1.6). In [3] Andrews' main focus is the application of the conjugate pair (1.5)--(1.6) for both $(\gamma_{n_1}(1),\delta_{n_1}(1))$ and $(\gamma_{n_2}(1),\delta_{n_2}(1))$ in (1.7) with $\rho_1,$ $\rho_2\rightarrow\infty.$
In [10] $q$-series were discovered related to both positive definite quadratic forms and indefinite ternary quadratic forms using only the 
$1$-fold Bailey lemma. A nice consequence of [10] is that it is clear that identities that arise naturally from inserting $2$-fold Bailey pairs into the $2$-fold Bailey lemma
may also be obtained using Bailey pairs with the $1$-fold Bailey lemma. One particularly nice example from that study is the new expansion
\begin{equation}(q)_{\infty}^3=\sum_{N\in\mathbb{Z}}(-1)^{N}q^{N(N-1)/2}\sum_{\substack{ n=0 \\ 2|j|\le n}}^{\infty}(-1)^{n+j}q^{n(n+1)/2-j(3j-1)/2+jN}.\end{equation}
In fact one may summarize the simple method from [10] in the following proposition.
\\*
{\bf Proposition 1.} \it If $(\gamma_n, \delta_n)$ is a conjugate Bailey pair relative to $a=1,$ then $(\alpha_N, \beta_N)$ is a Bailey pair relative to
$a=q$ where
\begin{equation}\alpha_N=\frac{1-q^{2N+1}}{1-q}
\sum_{j\in\mathbb{Z}} \gamma_j (-1)^{N+j} q^{\binom{j+N}{2}+j},\end{equation}
\begin{equation}\beta_N=\frac{\delta_Nq^{-N}}{(q)_{2N}}.\end{equation}

\rm

\par Given [10], one might wonder if there is a direct relation between a $1$-fold Bailey pair given a $2$-fold Bailey pair. As it turns out, there are many such relations and, in general, we have the following relationship between the $1$-fold and the $(s+1)$--fold Bailey pair.
\\*
\\*
{\bf Lemma 2.} \it Let $(A_{n_1,\dots,n_s,n},B_{n_1,\dots,n_s,n})$ form an $(s+1)$-fold Bailey pair relative to $(a_{n_1},a_{n_2},\dots,a_{n_{s}},a),$
and let $(\gamma^{(k)},\delta^{(k)})$ for $k=1,\dots,s$ be $s$ conjugate Bailey pairs relative to $a_{n_k}.$
Then $(A_n,B_n)$ defined by
\begin{equation}
\begin{aligned}
A_n&:=\sum_{n_1,\dots,n_s\in\mathbb{Z}} \gamma^{(1)}_{n_1}\cdots\gamma^{(s)}_{n_s} A_{n_1,\dots,n_s,n} \\
B_n&:=\sum_{n_1,\dots,n_s=0}^{\infty} \delta^{(1)}_{n_1}\cdots\delta^{(s)}_{n_s} B_{n_1,\dots,n_s,n}
\end{aligned}
\end{equation}
forms a $1$-fold Bailey pair relative to $a.$
\\*
\begin{proof} We write out the $s=1$ case. We have,
\begin{equation}\begin{aligned}&\sum_{n_1=0}^{\infty} \delta^{(1)}_{n_1}(a_1) B_{n_1,n}(a_1,a)\\ &=\sum_{n_1=0}^{\infty} \delta^{(1)}_{n_1}(a_1) \sum_{r_1=-n_1}^{n_1}\sum_{r_2=-n}^{n} \frac{A_{r_1,r_2}}{(a_1q)_{n_1+r_1} (q)_{n_1-r_1}(aq)_{n+r_2} (q)_{n-r_2}} \\ 
&=\sum_{r_2=-n}^{n} \frac{1}{(aq)_{n+r_2} (q)_{n-r_2}} \sum_{n_1=0}^{\infty} \delta^{(1)}_{n_1}(a_1)\sum_{r_1=-n_1}^{n_1}\frac{A_{r_1,r_2}}{(a_1q)_{n_1+r_1} (q)_{n_1-r_1}}\\ 
&=\sum_{r_2=-n}^{n} \frac{1}{(aq)_{n+r_2} (q)_{n-r_2}} \sum_{r_1=-\infty}^{\infty}A_{r_1,r_2}\sum_{n_1=|r_1|}^{\infty}\frac{\delta^{(1)}_{n_1}(a_1)}{(a_1q)_{n_1+r_1} (q)_{n_1-r_1}}\\ 
&=\sum_{r_2=-n}^{n} \frac{1}{(aq)_{n+r_2} (q)_{n-r_2}} \sum_{r_1=-\infty}^{\infty}A_{r_1,r_2}\gamma^{(1)}_{r_1}(a_1)\\
&=\sum_{r_2=-n}^{n} \frac{A_{r_2}(a)}{(aq)_{n+r_2} (q)_{n-r_2}}=B_n(a). \\
\end{aligned}\end{equation}
Since these steps can be repeated separately for each $k,$ the more general case follows.
\end{proof}

\rm
Lemma 2 appears to be new, and allows one to prove all of Andrews' pentagonal number theorem identities in [3] using only the $1$-fold Bailey lemma, by choosing $\gamma^{(k)}_{n_k}=q^{n_k^2}/(q)_{\infty},$ and $\delta^{(k)}_{n_k}=q^{n_k^2},$ for every $k$ between $1$ and $s.$

As an application of Lemma 2 we will show that Slater's well known
Bailey pair (A5), see [13, p.~463], follows from the ``diagonal''
$2$-fold Bailey pair [3, Eqs.~(4.1) and (4.4)]
\begin{equation}
\begin{aligned}
A_{n_1,n_2} &= (-1)^{n_1+n_2} q^{\binom{n_1+n_2}{2}} \\
B_{n_1,n_2} &=
\begin{cases} \displaystyle
\frac{1}{(q)_{2n_1}} & \text{if $n_1=n_2$} \\[2mm]
0 & \text{otherwise}
\end{cases}
\end{aligned}
\end{equation}
relative to $1$.
First we obtain Slater's $B_n$:
\begin{equation}\begin{aligned}
B_n &= \sum_{k=0}^{\infty} B_{k,n} q^{k^2} \\
 &= \sum_{k=0}^{\infty} \chi(k=n) \frac{q^{k^2}}{(q)_{2n}} \\
 &= \frac{q^{n^2}}{(q)_{2n}}.
\end{aligned} \end{equation}
For $A_n$ we have to work a little bit harder:

\begin{equation}\begin{aligned}
A_n &= \frac{1}{(q;q)_{\infty}} \sum_{k=-\infty}^{\infty} A_{k,n} q^{k^2} \\
&= \frac{(-1)^nq^{n(n-1)/2}}{(q;q)_{\infty}}\sum_{k=-\infty}^{\infty}(-1)^kq^{k(3k-1)/2+nj}.
\end{aligned} \end{equation}

By the Jacobi triple product identity [8, Eq.~(1.6.1)]
\[
\sum_{j=-\infty}^{\infty} (-a)^j q^{\binom{j}{2}}
=(q)_{\infty}(a)_{\infty}(q/a)_{\infty}
\]
this yields
\[
A_n = \frac{(-1)^n q^{n(n-1)/2}
 (q^3,q^3)_{\infty}(q^{n+1};q^3)_{\infty}(q^{2-n};q^3)_{\infty}}
 {(q;q)_{\infty}} .
\]
considering the three congruence classes of $n$ modulo $3$ this finally
simplifies to
\begin{equation}\alpha_n=\begin{cases}
1 & \text{if $n=0$} \\
q^{3k^2+k}+q^{3k^2-k} &\text{if $n=3k$} \\
-q^{3k^2  \pm k} &\text{if $n=3k \pm 1$}. \\
\end{cases}\end{equation}
\\*

Now given our computations, we can exploit the uniqueness of Bailey pairs to obtain a $2$-fold Bailey pair from Slater's $A(3)$ Bailey pair. We show this by first observing that the $\alpha_n$ corresponding to Slaters [13, A(3)] is the same as $\alpha_n$ in (1.16) but with $q$ replaced by $q^2,$
\begin{equation}\alpha_n=\begin{cases}
1 & \text{if $n=0$} \\
q^{6k^2+2k}+q^{6k^2-2k} &\text{if $n=3k$} \\
-q^{6k^2  \pm 2k} &\text{if $n=3k \pm 1$}. \\
\end{cases}\end{equation} 
That is, suppose $$B_n=\frac{q^n}{(q)_{2n}},$$ and in Lemma 2, set $s=1,$ and $\gamma^{(1)}_{n_1}=q^{n_1^2}/(q)_{\infty},$ and $\delta^{(1)}_{n_1}=q^{n_1^2}.$ The uniqueness of Bailey pairs with (1.17) and (1.15) tells us we must have
 $$\begin{aligned}A_n&=\frac{(-1)^nq^{n(n-1)}}{(q^2;q^2)_{\infty}}\sum_{j\in\mathbb{Z}}(-1)^jq^{j(3j-1)+2nj}\\
 &=\frac{1}{(q)_{\infty}}\sum_{r=-\infty}^{\infty}q^{r^2}A_{r,n},\end{aligned}$$
 and 
 $$\begin{aligned}B_n&=\frac{q^{n}}{(q)_{2n}}\\
 &=\sum_{r=0}^{\infty}q^{r^2}B_{r,n}.\end{aligned}$$
 
 Therefore, the $A_{r, n}$ and $B_{r, n}$ are forced once $(\gamma^{(1)}_{n_1},\delta^{(1)}_{n_1})$ is chosen, and we have proven that we must have the symmetric $2$-fold Bailey pair
\begin{equation}A_{n_1,n_2}=\frac{(-1)^{n_1+n_2}}{(-q)_{\infty}}q^{4\binom{n_1}{2}+n_1+2\binom{n_2}{2}+2n_1n_2},\end{equation}
\begin{equation}B_{n_1,n_2}=\begin{cases} 0,& \text {if } n_2\neq n_1,\\ \displaystyle\frac{q^{n_1-n_1^2}}{(q)_{2n_1}}, & \text{if } n_1= n_2.\end{cases}\end{equation}

This pair gives us the identity due to L.J. Rogers [12, pg.~332, Eq.~(13)] (from the $\rho_1,\rho_2\rightarrow\infty$ case of (1.5)--(1.6))
\begin{equation}\begin{aligned}&\sum_{n=0}^{\infty}\frac{q^{n^2+n}}{(q)_{2n}}=\frac{1}{(q)_{\infty}(q^2;q^2)_{\infty}}\sum_{n\in\mathbb{Z}}(-1)^nq^{n(2n-1)}\sum_{j\in\mathbb{Z}}(-1)^jq^{j(3j-1)+2nj}\\ &=\frac{1}{(q^3,q^4,q^5,q^6,q^7;q^{10})_{\infty}(q^2,q^{18};q^{20})_{\infty}}.\end{aligned}\end{equation}
It is important to note that the literature has a large volume of $1$-fold Bailey pairs, and so our argument in obtaining (1.18)--(1.19) is a more natural and potent strategy in obtaining further $2$-dimensional identities. 
\par Andrews used (1.13) to obtain a nice two-dimensional pentagonal number theorem identity [3, Theorem 2]. For more $2$-fold and $3$-fold Bailey pairs see Berkovich [7].
\par Tactically speaking, the pair (1.14)--(1.15) presently has more utility, as there are more known conjugate Bailey pairs for the $1$-fold Bailey lemma (e.g. [4]). However, (1.13) would appear to encompass a larger pool of identities overall, as (1.14)--(1.15) is obtained from the limiting case $\rho_1,\rho_2\rightarrow\infty$ of (1.5)--(1.6) with (1.13).

\section{The $\spt(n)$ function of Andrews}
In [11], we encountered the double sum
$$\sum_{n_1,n_2=1}^{\infty}\frac{q^{n_1+n_2}}{(1-q^{n_1})^2(q^{n_1+1})_{\infty}(1-q^{n_2})^2\dots(1-q^{n_1+n_2})},$$
and asked if the sum over $n_2$ had any origin from the $1$-fold Bailey lemma. The proof relied on a $2$-fold Bailey pair from [9]. It was also suggested there was a new generalized form of
Andrews' relation [5] $\spt(n)=np(n)-\frac{1}{2}N_2(n).$ Here $\spt(n)$ is the total number of appearances of the smallest parts of all the partitions of $n,$ $p(n)$ is the classical
unrestricted partition function, and $N_2(n)$ is the second Atkin-Garvan moment (see (2.6) and [5] for the generating function).
\\*
{\bf Lemma 3.} \it For $n$ and $M$ non-negative integers, $(\alpha_n, \beta_n)$ forms a Bailey pair relative to $a=1,$ where
\begin{equation}\alpha_n=\frac{(q)_{M}^2(-1)^n(1+q^n)q^{n(3n-1)/2}}{(q)_{M-n}(q)_{M+n}},\end{equation}
for $1\leq n\leq M,$ $\alpha_n=0$ if $n>M,$
$$\alpha_0=1,$$
and
\begin{equation}\beta_n=\frac{(q)_{M}}{(q)_n(q)_{n+M}}.\end{equation}

\begin{proof} Using the inverse relation of a Bailey pair [2] (or [15, Eq.~(2.4)]), 
\begin{equation}\alpha_n=\frac{(1-aq^{2n})(a)_n(-1)^nq^{n(n-1)/2}}{(1-a)(q)_n}\sum_{k=0}^{n}(q^{-n})_k(aq^n)_kq^k\beta_k,\end{equation}
we choose our $\beta_n$ to be (2.2), insert into (2.3), and write
$$(-1)^N(1+q^N)q^{\binom{N}{2}}\sum_{j=0}^{N}\frac{(q^N)_j(q^{-N})_jq^j}{(q)_j(q)_{j+M}}=\frac{(q)_{M}(-1)^N(1+q^N)q^{N(3N-1)/2}}{(q)_{M-N}(q)_{M+N}}=\alpha_N,$$
for $N>0,$ and 
$$\alpha_0=\frac{1}{(q)_M}.$$
This follows from the $q$--Chu--Vandermonde theorem [8, Eq.~(II.6)] with $(a,c,n)\mapsto (q^N,q^{M+1},N),$ because
$$(q)_M\sum_{j=0}^{N}\frac{(q^N)_j(q^{-N})_jq^j}{(q)_j(q)_{j+M}}=\frac{(q^{M-N+1})_Nq^{N^2}}{(q^{M+1})_N}=\frac{(q)_M(q)_Mq^{N^2}}{(q)_{M-N}(q)_{M+N}}.$$
Finally, to respect the convention that $\alpha_0=1$, we multiply through by $(q)_{M}.$ 
\end{proof}

{\bf Corollary 4.} We have, for each natural number $M,$
\begin{equation}\sum_{n=1}^{\infty}\frac{q^n}{(1-q^n)^2(1-q^{n+1})\cdots(1-q^{n+M})} \end{equation}
$$=\frac{1}{(q)_M}\sum_{n=1}^{\infty}\frac{nq^n}{1-q^n}+(q)_{M}\sum_{n=1}^{M}\frac{(-1)^n(1+q^n)q^{n(3n+1)/2}}{(q)_{M-n}(q)_{M+n}(1-q^n)^2}.$$

\rm
\begin{proof} Using the conjugate pair (1.5)--(1.6), differentiating with respect to $\rho_1,\rho_2$ and then putting $\rho_1,\rho_2=1,$ we obtain
\begin{equation}\sum_{n=1}^{\infty}(q;q)_{n-1}^2\beta_nq^n=\alpha_0\sum_{n=1}^{\infty}\frac{nq^n}{1-q^n}+\sum_{n=1}^{\infty}\frac{\alpha_nq^n}{(1-q^n)^2}.\end{equation}
Applying the Bailey pair contained in Lemma 3 to (2.5) now gives the theorem. \end{proof}

\par The $q$-series on the left side of (2.4) appeared in [11], and is the generating function for $\spt_M^{*}(n),$ the total number of appearances of the smallest parts
of the number of partitions of $n$ where parts greater than the smallest plus $M$ do not occur. The first sum on the right side of (2.4) may be interpreted as $\sum_{k=0}^{n}\sigma_1(k)p_M(n-k),$ where $\sigma_1(n)=\sum_{d|n}d,$ and $p_M(n)$ is the number of partitions of $n$ into parts $\le M.$ The limiting case $M\rightarrow\infty$ is Euler's well known formula $np(n)=\sum_{k=0}^{n}\sigma_1(k)p(n-k),$ which is also an observation used by Andrews to obtain his $\spt(n)$ identity [5]. By Tannery's theorem [14, pg.~292] and [5, Eq.~(3.4)], it can be seen that the limit of the second sum on the right side of Corollary 4 is
\begin{equation}\frac{1}{2}\sum_{n=1}^{\infty}N_2(n)q^{n}=-\frac{1}{(q)_{\infty}}\sum_{n=1}^{\infty}\frac{(-1)^nq^{n(3n+1)/2}(1+q^n)}{(1-q^{n})^2}.\end{equation}
The second sum on the right hand side of Corollary 4 is more complicated, and is worthy of a separate study, as any information on $\spt^{*}_M(n)$ is important to better understand $\spt(n).$ The case $M\rightarrow\infty$ of Corollary 4 can now be seen as $\spt(n)=np(n)-\frac{1}{2}N_2(n).$ While Corollary 4 is important in its own right, it also implies the following Bailey pair.
\\*
\\*
{\bf Lemma 5.} \it For $n$ and $M$ non-negative integers, $(\alpha_M, \beta_M)$ forms a Bailey pair relative to $a=1,$ where
\begin{equation}\alpha_M=\frac{(-1)^M(1+q^M)q^{M(3M+1)/2}}{(1-q^M)^2},\end{equation}
for $M>0,$ and 
$$\alpha_0=\sum_{n=1}^{\infty}\frac{nq^n}{1-q^n},$$
\begin{equation}\beta_M=\frac{1}{(q)_{M}}\sum_{n=1}^{\infty}\frac{q^n}{(1-q^n)^2(1-q^{n+1})\cdots(1-q^{n+M})}.\end{equation}
\\* \rm  \par 
The point of this section is that no discussion of $\spt_{M}^{*}(n)$ (or Corollary 4) arose until studying some identities using the $2$-fold Bailey lemma.

1390 Bumps River Rd. \\*
Centerville, MA
02632 \\*
USA \\*
E-mail: alexpatk@hotmail.com, alexepatkowski@gmail.com
\end{document}